\documentclass[12pt]{article}
\usepackage{etex}
\usepackage{setspace}
\usepackage[british]{babel}
\usepackage{microtype}
\usepackage{amssymb}
\usepackage{appendix} 
\usepackage{amsmath}
\usepackage{graphs}
\usepackage{times}
\usepackage{anysize} 
\usepackage[Lenny]{fncychap}
\usepackage{lettrine}
\usepackage{bigstrut,array,float, fancyvrb,marvosym,verbatim,combgames}
\marginsize{25mm}{25mm}{40mm}{25mm}

\usepackage{psgo,array,url}

\newtheorem{definition}{Definition}[section]

\newtheorem{theorem}{Theorem}[section]
\newtheorem{conjecture}{Conjecture}[section]

\newcommand{\qed}{\hfill \bf q.e.d.} 
\newcommand{\lef}{\mathcal{L}}
\newcommand{\ri}{\mathcal{R}}
\newcommand{\n}{\mathcal{N}}
\newcommand{\pre}{\mathcal{P}}
\newcommand{\ti}{\mathcal{T}} 
\newcommand{\gfr}{G_F^{SR}}
\newcommand{\gfl}{G_F^{SL}}
\newcommand{\plus}{+_{\ell}}

\newcommand{\G}{\mathcal{G}}
\newenvironment{proof}{{\noindent\it Proof}: }{\qed}

\makeatletter
{\catcode`\|=\z@ \catcode`\\=12 |gdef|bslash{\}}
\makeatother

\onehalfspacing

\setgounit{0.4cm}

\begin{document}
\title{Scoring Play Combinatorial Games Under Different Operators}
\author{Fraser Stewart\\
\texttt{fraseridstewart@gmail.com}}
\date{}

\maketitle{}

\centerline{\bf Abstract}
Scoring play games were first studied by Fraser Stewart for his PhD thesis \cite{FSP}.  He showed that under the disjunctive sum, scoring play games are partially ordered, but do not have the same ``nice'' structure of normal play games.  In this paper I will be considering scoring play games under three different operators given by John Conway \cite{ONAG} and William Stromquist and David Ullman \cite{SU}, namely the conjunctive sum, selective sum and sequential join.

\section{Introduction}

Until very recently scoring play games have not received the kind of treatment or analysis that normal and mis\`ere play games have.  The general definition of a scoring play game is given below, for further reading on the general structure of scoring play games see \cite{FS} and \cite{FSP}.

In this paper we will be examining scoring play games under three different operators; the conjunctive sum, where the players must move on all available components on their turn; the selective sum, where the players can pick any components they wish to move on on their turn; and finally the sequential join, where the components are given a pre-arranged order and the players must play on them in that order.  These operators were first defined by John Conway in On Numbers and Games \cite{ONAG}, and William Stromquist and David Ullman \cite{SU}.

We will also be looking at the Sprague-Grundy values of scoring play octal games under these three different operators.  We will give evidence and conjecture that the period of the scoring play Sprague-Grundy function is eventually periodic, and has the same period for finite octal games as the disjunctive sum, for all three operators.

\subsection{Scoring Play Theory}

Intuitively a scoring play game is one that has the following three properties;

\begin{enumerate}
\item{The rules of the game clearly define what points are and how players either gain or lose them.}
\item{When the game ends the player with the most points wins.}
\item{For any two games $G$ and $H$, $a$ points in $G$ are equal to $a$ points in $H$, where $a\in \mathbb{R}$.}
\item{At any stage in a game $G$ if Left has $L$ points and Right has $R$ points, then the score of $G$ is $L-R$, where $L,R\in \mathbb{R}$.}
\end{enumerate}

Mathematically the definition is given as follows \cite{FS};

\begin{definition} A scoring play game $G=\{G^L|G^S|G^R\}$, where $G^L$ and $G^R$ are sets of games and $G^S\in\mathbb{R}$, the base case for the recursion is any game $G$ where $G^L=G^R=\emptyset$.\\

\noindent $G^L=\{\hbox{All games that Left can move to from } G\}$\\
\noindent $G^R=\{\hbox{All games that Right can move to from } G\}$,\\

and for all $G$ there is an $S=(P,Q)$ where $P$ and $Q$ are the number of points that Left and Right have on $G$ respectively.  Then $G^S=P-Q$, and for all $g^L\in G^L$, $g^R\in G^R$, there is a $p^L,p^R\in\mathbb{R}$ such that $g^{LS}=G^S+p^L$ and $g^{RS}=G^S+p^R$.

$\gfl$ and $\gfr$ are called the final scores of $G$ and are the largest scores that Left and Right can achieve when $G$ ends, moving first respectively, if both players play their optimal strategy on $G$.

\end{definition}

For scoring play the disjunctive sum needs to be defined a little differently, because in scoring games when we combine them together we have to sum the games and the scores separately.  For this reason we will be using two symbols $\plus$ and $+$.  The $\ell$ in the subscript stands for ``long rule'', this comes from \cite{ONAG}, and means that the game ends when a player cannot move on any component on his turn.

\begin{definition}  The disjunctive sum is defined as follows:
$$G\plus H=\{G^L\plus H,G\plus H^L|G^S+H^S|G^R\plus H,G\plus H^R\},$$
\noindent where $G^S+H^S$ is the normal addition of two real numbers.
\end{definition}

The outcome classes also need to be redefined to take into account the fact that a game can end with a tied score.  So we have the following two definitions.

\begin{definition}
\item{$L_>=\{G|\gfl>0\}$, $L_<=\{G|\gfl<0\}$, $L_= =\{G|\gfl=0\}$.}
\item{$R_>=\{G|\gfr>0\}$, $R_<=\{G|\gfr<0\}$, $R_= =\{G|\gfr=0\}$.}
\item{$L_\geq=L_>\cup L_=$, $L_\leq = L_<\cup L_=$.}
\item{$R_\geq=R_>\cup R_=$, $L_\leq = R_<\cup R_=$.}
\end{definition}

\begin{definition}
The outcome classes of scoring games are defined as follows:

\begin{itemize}
\item{$\lef=(L_>\cap R_>)\cup(L_>\cap R_=)\cup(L_=\cap R_>)$}
\item{$\ri=(L_<\cap R_<)\cup(L_<\cap R_=)\cup(L_=\cap R_<)$}
\item{$\n=L_>\cap R_<$}
\item{$\pre=L_<\cap R_>$}
\item{$\ti=L_=\cap R_=$}
\end{itemize}
\end{definition}

We will also be using two conventions throughout this paper.  The first is that the initial score of a game will be $0$ unless stated otherwise.  The second is that for a game $G$ if $G^L=G^R=\emptyset$, then we will write $G$ as $G^S$ rather than $\{.|G^S|.\}$.  For example the game $G = \{\{.|0|.\}|1|\{.|2|.\}\}$, will be written as $\{0|1|2\}$. The game $\{.|n|.\}$, will be written as $n$, and so on. This is simply for convenience and ease of reading.

\subsection{Impartial Games}

The definition of an impartial scoring play game is less intuitive than for normal and mis\`ere play games.  The reason for this is because we have to take into account the score, for example, consider the game $G=\{4|3|2\}$.  On the surface the game does not appear to fall into the category of an impartial game, since Left wins moving first or second, however this game is impartial since both players move and gain a single point, i.e. they both have the same options.

So we will use the following definition for an impartial game;

\begin{definition}
A scoring game $G$ is impartial if it satisfies the following;

\begin{enumerate}
\item{$G^L=\emptyset$ if and only if $G^R=\emptyset$.}
\item{If $G^L\neq \emptyset$ then for all $g^L\in G^L$ there is a $g^R\in G^R$ such that $g^L\plus -G^S=-(g^R\plus -G^S)$.}
\end{enumerate}
\end{definition} 

We will also be looking at octal games in this paper, and for scoring play games we use the following definition of an octal game.

\begin{definition}
A scoring play octal game $O=(n_1n_2\dots n_k, p_1p_2\dots p_k)$, is a set of rules for playing nim where if a player removes $i$ beans from a heap of size $n$ he gets $p_i$ points, $p_i\in \mathbb{R}$, and he must leave $a,b,c\dots$ or $j$ heaps, where $n_i=2^a+2^b+2^c+\dots+2^j$.
\end{definition}

By convention we will also say that $n\in O$ means that the nim heap $n$ is played under the rule set $O$.  In \cite{FSP} and \cite{FSI}, the following definition and conjecture were given.  

\begin{definition}
Let $n\in O=(t_1t_2\dots t_f, p_1p_2\dots p_t)$ and $m\in P=(s_1s_2\dots s_e, q_1q_2\dots q_t)$;

\begin{itemize}
\item{$\G_s(0)=0$.}
\item{$\G_s(n)=\max_{k,i}\{p_k-\G_s(n_1\plus n_2\plus \dots\plus n_{i})\}$, where $n_1+n_2+\dots +n_{i}=n-k$, $t_k=\Sigma_{i\in S_k}2^i$.}
\item{$\G_s(n\plus m)=\max_{k,i,l,j}\{p_k-\G_s(n_1\plus n_2\plus \dots\plus n_i\plus m),q_l-\G_s(n\plus m_1\plus m_2\plus\dots\plus m_j)\}$, where $n_1+n_2+\dots +n_i=n-k$, $t_k=\Sigma_{i\in S_k}2^i$, $m_1+m_2+\dots m_j=m-l$ and $s_l=\Sigma_{j\in R_l}2^j$.}
\end{itemize}
\end{definition}

\begin{conjecture}\label{period}
Let $O=(n_1n_2\dots n_t, p_1p_2\dots p_t)$ and $P=(m_1m_2\dots m_l, q_1q_2\dots q_l)$ be two finite taking-no-breaking octal games such that, there is at least one $n_s\neq 0$ or $1$, and if $n_i$ and $m_j=1,2$ or $3$ then $p_i=i$ and $q_j=j$, and $p_i=q_j=0$, otherwise, then for all $m$ there exists an $N$ such that;

$$\G_s(n+2k\plus m)=\G_s(n\plus m)$$

\noindent for all $n\geq N$and $k$ is the largest entry in $O$ such that $n_k\neq 0,1$.
\end{conjecture}

There was a lot of evidence given that the conjecture is true.  In this thesis we will be examining the same function under all three operators.  I also conjecture that if conjecture \ref{period} is true then the function settles down to the same period for all 3 operators.  If this is true I think it would be a very interesting result, since you would expect changing the operator would change the period of the function, but in this case it appears that this probably does not happen.

\section{The Conjunctive Sum}

The first operator that we will be looking at is the conjunctive sum.  Under this operator, players must move on all components on their turn.  Mathematically it is defined as follows;

\begin{definition} The conjunctive sum is:

$$G\bigtriangleup H=\{G^L\bigtriangleup H^L|G^S+H^S|G^R\bigtriangleup H^R\}$$

\noindent where $G^S + H^S$ is the normal addition of two real numbers.
\end{definition}

\begin{theorem}
If $G\not\cong 0$ then $G\neq 0$.
\end{theorem}

\begin{proof}  First consider the game $G^L=G^R=\emptyset$, then clearly if $G^S\neq 0$ then $G\neq 0$.

Next consider the case where $G^L\neq \emptyset$, since the case $G^R\neq \emptyset$ follows by symmetry.  Let $P=\{.|a|b\}$, where $a=P^{SL}_F>0$.  Since $G$ is a combinatorial game, this means that the game tree has finite depth and finite width, so we can choose $b<0$ such that $|b|$ is greater than any number on $G$.  On Left's first turn he must move to $G^L\bigtriangleup P$, regardless of whether Right can play $G$ or not, he will have to move on $P$ on his next turn.

Thus $(G\bigtriangleup P)^{SL}_F<0$, and therefore $G\bigtriangleup P\not\approx P$, and the theorem is proven.  
\end{proof}

\begin{theorem}
For any outcome classes $\mathcal{X}$, $\mathcal{Y}$ and $\mathcal{Z}$, there is a game $G\in \mathcal{X}$ and $H\in \mathcal{Y}$ such that $G\bigtriangleup H\in \mathcal{Z}$.
\end{theorem}

\begin{proof}  To prove this consider the following game $G\bigtriangleup H$, where $G=\{\{.|b|\{.|c|\{e|d|f\}\}\}|a|.\}$ and $H=\{.|g|\{\{\{k|j|.\}|i|.\}|h|.\}\}$, as shown in figure \ref{cout}.

\begin{figure}[htb]
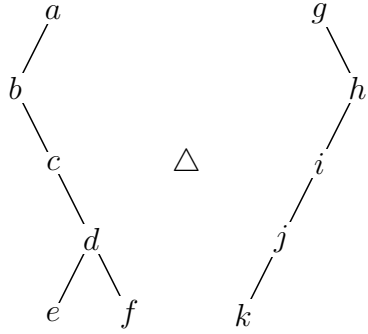

\begin{center}
\begin{graph}(4.5,4)
\roundnode{1}(0.5,4)\roundnode{2}(0,3)\roundnode{3}(0.5,2)\roundnode{4}(1,1)\roundnode{5}(1.5,0)
\roundnode{6}(0.5,0)\roundnode{7}(3,0)\roundnode{8}(3.5,1)\roundnode{9}(4,2)\roundnode{10}(4.5,3)
\roundnode{11}(4,4)

\edge{1}{2}\edge{2}{3}\edge{3}{4}\edge{4}{5}\edge{4}{6}
\edge{7}{8}\edge{8}{9}\edge{9}{10}\edge{10}{11}

\freetext(2.25,2){$\bigtriangleup$}

\nodetext{1}{$a$}\nodetext{2}{$b$}\nodetext{3}{$c$}\nodetext{4}{$d$}\nodetext{5}{$f$}
\nodetext{6}{$e$}\nodetext{7}{$k$}\nodetext{8}{$j$}\nodetext{9}{$i$}\nodetext{10}{$h$}
\nodetext{11}{$g$}
\end{graph}
\end{center}
\caption{$\{\{.|b|\{.|c|\{e|d|f\}\}\}|a|.\}\bigtriangleup \{.|g|\{\{\{k|j|.\}|i|.\}|h|.\}\}$}\label{cout}
\end{figure}

In these games $G^{SL}_F=c$ and $G^{SR}_F=a$, $H^{SL}_F=g$ and $H^{SR}_F=i$, however $(G\bigtriangleup H)^{SL}_F=e+j$ and $(G\bigtriangleup H)^{SR}_F=e+k$.  Since the outcome classes of $G$ and $H$ depend on $a$, $c$, $g$ and $i$, and the outcome class of $G\bigtriangleup H$ depends on $e+j$ and $e+k$, then clearly we can choose $a$, $c$, $g$, $i$, $e$, $j$ and $k$, so that $G$ and $H$ can be in any outcome class and $G\bigtriangleup H$ can be in any outcome class and the theorem is proven.
\end{proof}

\subsection{Impartial Games}

\begin{theorem}
Impartial games form an abelian group under the conjunctive sum.
\end{theorem}

\begin{proof}  To prove this we only need to show that there is an identity set $I$ that contains more than one element, and that for any impartial game $G$, there is a $G^{-1}$ such that $G\bigtriangleup G^{-1}\in I$.

Let $I=\{G|G \hbox{ is impartial and }G\in \ti\}$, then we wish to show that for all $G\in I$, $G\bigtriangleup P\approx P$ for all impartial games $P$.  There are three cases to consider, since the remaining follow by symmetry, $P\in L_>$, $P\in L_<$ or $P\in L_=$.

So first let $P\in L_>$, and consider the game $G\bigtriangleup P$.  Since Left can achieve a score of $0$ on $G$, then all Left has to do is play his winning strategy on $P$, and $G\bigtriangleup P\in L_>$.

Next let $P\in L_<$, and consider the game $G\bigtriangleup P$.  $G\in L_=$, and since both $G$ and $P$ are impartial, neither player can change the parity of either game, since they must both play both games on every turn.     So all Right has to do is play his winning strategy on $P$ and $G\bigtriangleup P\in L_<$.

Finally let $P\in L_=$, and consider the game $G\bigtriangleup P$.  If both players always make their best moves on $G$ and $P$ then the final score of $G\bigtriangleup P$ will be 0, since $G\in \ti$ and $P\in L_=$.  Since $G\in \ti$ and $P\in L_=$, this implies that if Left chooses a different move other than his best move either $G$ or $P$, then the final score will be $\leq 0$, and similarly for Right.  This means that as long as Right keeps playing his best strategy, if Left chooses anything else Right can potentially win and similarly for Left.  In other words the best thing for both players to do is to play their best strategy on both $G$ and $P$ and the final score will be a tie, i.e.  $G\bigtriangleup P\in L_=$.

The cases for $R_>$, $R_<$ and $R_=$ follow by symmetry.

For the inverse of a game $G$, where $G^{SL}_F=n$ and $G^{SR}_F=p$, we let $H$ be a game where $H^{SL}_F=-n$ and $H^{SR}_F=-p$.  Note that $G\bigtriangleup H\in I$ if and only if $G\bigtriangleup H\in \ti$.

So consider the game $G\bigtriangleup H$ with Left moving first, since the case where Right moves first follows by symmetry.  If $G^L=\emptyset$, then this implies that $G^R=\emptyset$ since $G$ is impartial, which implies that $G^{SL}_F=G^{SR}_F=n$, so for the inverse let $H$ be a game such that $H^{SL}_F=H^{SR}_F=-n$.  However since $H$ is impartial the only game that satisfies that condition is the game $H=\{.|-n|.\}$, which is clearly the inverse of $G$.

If $G^L\neq \emptyset$, and Left and Right make their best move at every stage on both $G$ and $H$, then the final score of $G\bigtriangleup H$ will be $G^{SL}_F+H^{SL}_F=n-n=0$.  Using the same argument as the identity proof if Left or Right try a different strategy then the final score will be either $\leq 0$ or $\geq 0$ respectively, therefore $(G\bigtriangleup H)^{SL}_F=(G\bigtriangleup H)^{SR}_F=0$ and $H$ is the inverse of $G$.

It is clear that the set is closed, since if $G$ and $H$ are impartial then $G\bigtriangleup H$ must also be impartial.  It is also clear that we have commutativity and associativity, since we must play on every component on every turn, then the order of the components is irrelevant.
\end{proof}

\subsection{Sprague-Grundy Theory}

First we define the following;

\begin{definition}
Let $n\in O=(t_1t_2\dots t_f, p_1,\dots p_f)$ and $m\in P=(s_1s_2\dots s_e, q_1\dots q_e)$;

\begin{itemize}
\item{$\G_s(0)=0$.}
\item{$\G_s(n)=\max_{k,i}\{p_k-\G_s(n_1\bigtriangleup n_2\bigtriangleup \dots\bigtriangleup n_{i})\}$, where $n_1+n_2+\dots +n_{i}=n-k$ and $t_k=\Sigma_{i\in S_k}2^i$.}
\item{$\G_s(n\plus m)=\max_{k,i,l,j}\{p_k+q_l-\G_s(n_1\bigtriangleup n_2\bigtriangleup \dots\bigtriangleup n_i\bigtriangleup m_1\bigtriangleup m_2\bigtriangleup\dots\bigtriangleup m_j)\}$, where $n_1+n_2+\dots +n_i=n-k$,$t_k=\Sigma_{i\in S_k}2^i$, $m_1+m_2+\dots m_j=m-l$ and $s_l=\Sigma_{j\in R_l}2^j$.}
\end{itemize}
\end{definition}

The fact that impartial games are a group mean that we can easily solve any octal game simply by knowing each heap's $\G_s(n)$ value.  So we have the following theorem;

\begin{theorem}
$$\G_s(n\bigtriangleup m)=\G_s(n)+\G_s(m)$$
\end{theorem}

\begin{proof}  We will prove this by induction.  The base case is trivial, since $\G_s(0\bigtriangleup 0)=\G_s(0)+\G_s(0)=0$.

So assume that the theorem holds for all values up to $\G_s(n\bigtriangleup m)$, and consider $\G_s(n+1\bigtriangleup m)$, since the case $\G_s(n\bigtriangleup m+1)$ follows by symmetry.  $\G_s(n+1\bigtriangleup m)=\max_{k,i,l,j}\{p_k+q_l-\G_s(n_1\bigtriangleup n_2\bigtriangleup \dots\bigtriangleup n_i\bigtriangleup m_1\bigtriangleup m_2\bigtriangleup\dots\bigtriangleup m_j)\}$, where $n_1+n_2+\dots +n_i=n+1-k$, $t_k=\Sigma_{i\in S_k}$, $m_1+m_2+\dots m_j=m-l$ and $s_l=\Sigma_{j\in R_l}$.  However each $n_i'<n+1$ and $m_j'<m$, therefore $\max_{k,i,l,j}\{p_k+q_l-\G_s(n_1\bigtriangleup n_2\bigtriangleup \dots\bigtriangleup n_i\bigtriangleup m_1\bigtriangleup m_2\bigtriangleup\dots\bigtriangleup m_j)\}=\max_{k,i,l,j}\{k+l-\G_s(n_1)+\G_s(n_2)+\dots+\G_s(n_i)+\G_s(m_1)+\G_s(m_2)+\dots+\G_s(m_j)\}$, by induction.  

Therefore $\max_{k,i,l,j}\{p_k+q_ll-\G_s(n_1)+\G_s(n_2)+\dots+\G_s(n_i)+\G_s(m_1)+\G_s(m_2)+\dots+\G_s(m_j)\}=\max_{k,i}\{p_k-G_s(n_1)+\G_s(n_2)+\dots+\G_s(n_i)\}+\max_{m,j}\{q_l-\G_s(m_1)+\G_s(m_2)+\dots+\G_s(m_j)\}=\max_{k,i}\{p_k-\G_s(n_1\bigtriangleup n_2\bigtriangleup\dots\bigtriangleup n_i)\}+\max_{m,j}\{q_l-\G_s(m_1\bigtriangleup m_2\bigtriangleup\dots\bigtriangleup m_j)\}=\G_s(n+1)+\G_s(m)$, and the proof is finished.
\end{proof}

\section{The Selective Sum}

The selective sum is a more general version of the disjunctive sum.  Rather than choosing a single component on each turn and playing that one only, the player can select any components he wishes to play and play those components on his turn.  It is defined as follows;

\begin{definition} The selective sum is:

$$G\triangledown H=\{G^L\triangledown H, G\triangledown H^L, G^L\triangledown H^L|G^S+H^S|G^R\triangledown H, G\triangledown H^R, G^R\triangledown H^R\}$$

\noindent where $G^S +H^S$ is the normal addition of two real numbers.
\end{definition}

\begin{theorem}
If $G\not\cong 0$ then $G\neq 0$.
\end{theorem}

\begin{proof}  The proof of this is very similar to the same theorem for the conjunctive sum.  First consider the game $G^L=G^R=\emptyset$, then clearly if $G^S\neq 0$ then $G\neq 0$.

Next consider the case where $G^L\neq \emptyset$, since the case $G^R\neq \emptyset$ follows by symmetry.  Let $P=\{.|a|b\}$, where $a=P^{SL}_F>0$.  Since $G$ is a combinatorial game, this means that the game tree has finite depth and finite width, we can choose $b$ to be more negative than any number on $G$.  On Left's first turn he must move to $G^L\triangledown P$, Right can then win by simply moving to $G^L\triangledown b$ on his turn, since the final score will be less than $0$, regardless of what Left does.

Thus $(G\triangledown P)^{SL}_F<0$, and therefore $G\triangledown P\not\approx P$, and the theorem is proven.     
\end{proof}

\begin{theorem}
For any outcome classes $\mathcal{X}$, $\mathcal{Y}$ and $\mathcal{Z}$, there is a game $G\in \mathcal{X}$ and $H\in \mathcal{Y}$ such that $G\bigtriangleup H\in \mathcal{Z}$.
\end{theorem}

\begin{proof}  To prove this consider the following game $G\triangledown H$, where $G= \{\{c|b|.\}|a|.\}$ and $H=\{.|d|\{.|e|\{.|f|g\}\}\}$, as shown in the following diagram.  

\begin{figure}[htb]
\begin{center}
\begin{graph}(3.5,3)

\roundnode{1}(0,1)\roundnode{2}(0.5,2)\roundnode{3}(1,3)\roundnode{4}(2,3)
\roundnode{5}(2.5,2)\roundnode{6}(3,1)\roundnode{7}(3.5,0)

\edge{1}{2}\edge{2}{3}\edge{4}{5}\edge{5}{6}\edge{6}{7}

\freetext(1.5,1.5){$\triangledown$}

\nodetext{1}{$a$}\nodetext{2}{$b$}\nodetext{3}{$c$}\nodetext{4}{$d$}
\nodetext{5}{$e$}\nodetext{6}{$f$}\nodetext{7}{$g$}

\end{graph}
\end{center}
\caption{$\{\{c|b|.\}|a|.\}\triangledown \{.|d|\{.|e|\{.|f|g\}\}\}$}
\end{figure}

In these games $G^{SL}_F=b$ and $G^{SR}_F=a$, $H^{SL}_F=d$ and $H^{SR}_F=e$, however $(G\triangledown H)^{SL}_F=c+f$ and $(G\bigtriangleup H)^{SR}_F=c+g$.  Since the outcome classes of $G$ and $H$ depend on $a$, $b$, $d$ and $e$, and the outcome class of $G\triangledown H$ depends on $c+f$ and $c+g$, then clearly we can choose $a$, $b$, $c$, $d$, $e$, $f$ and $g$, so that $G$ and $H$ can be in any outcome class and $G\triangledown H$ can be in any outcome class and the theorem is proven.

\end{proof}

\subsection{Impartial Games}

\begin{theorem}
Impartial games form a non-trivial monoid under the selective sum.
\end{theorem}

\begin{proof}  To prove that we have a non-trivial monoid we simply need to define an identity set that contains more than the game $\{.|0|.\}$.

First I will define a subset of the impartial games as follows;

$$I=\{i|G\plus i\approx G, \hbox{ for all impartial games }G\}$$  

Again, in order to show that we have a non-trivial monoid we have to show that $I$ contains more than one element.  So consider the following impartial game, $i=\{\{0|0|0\}|0|\{0|0|0\}\}$.

\begin{figure}[htb]
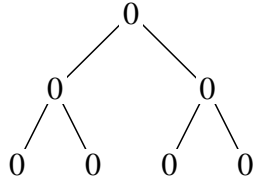

\begin{center}
\begin{graph}(3,2)

\roundnode{1}(0,0)\roundnode{2}(0.5,1)\roundnode{3}(1,0)\roundnode{4}(1.5,2)
\roundnode{5}(2,0)\roundnode{6}(2.5,1)\roundnode{7}(3,0)

\edge{1}{2}\edge{2}{3}\edge{4}{2}\edge{4}{6}\edge{5}{6}\edge{6}{7}

\nodetext{1}{0}\nodetext{2}{0}\nodetext{3}{0}\nodetext{4}{0}
\nodetext{5}{0}\nodetext{6}{0}\nodetext{7}{0}

\end{graph}
\end{center}
\caption{The game $\{\{0|0|0\}|0|\{0|0|0\}\}$}
\end{figure}

To show that $i\triangledown G\approx G$ for all impartial games $G$, there are 3 cases to consider $\gfl>0$, $\gfl<0$ and $\gfl=0$, since the cases for Right follow by symmetry.  First let $\gfl>0$, if Left has no move on $G$, then neither does Right, since $G$ is impartial, i.e. $G=G^S$, so they will play $i$ and the final score will still be $G^S$.

So let Left have a move on $G$, if Left chooses his best move on $G$,  then if Right plays $i$, then Left will respond in $i$ and Right must play $G$, which Left wins.  If Right tries to play on both $G$ and $i$, then either Right moves to a game where $G^L=G^R=\emptyset$, in which case Left moves on $i$ only and wins, or $G^L\neq \emptyset$, and Left also plays both $G$ and $i$ in order to maintain parity on $G$ and still wins.  Clearly if Right chooses to play $G$, then he will still lose, since Left also plays $G$ until it is finished and neither player can gain points on $i$.

Next let $\gfl<0$, this means that no matter what Left does, he will lose playing only $G$ on $G\triangledown i$, since Right will simply respond in $G$, until $G$ is finished, then they will play $i$, which does not change the final score of $G$.  Again if Left tries to change the parity of $G$, by playing $i$, Right will also play $i$, and it will be Left's turn to move on $G$ again.  If Left chooses to move on both $G$ and $i$, then as before Right will also move on $G$ and $i$ if $G^R\neq\emptyset$, and $i$ if $G^R=\emptyset$, but will win either way.

Finally let $\gfl=0$.  This means that Left's best move will be a move that eventually ties $G$.  So consider the game $G\bigtriangleup i$, Left's best move will be to move either on $G$ or $G$ and $i$, if Left moves on $i$ then this will give Right an opportunity to move first on $G$ and potentially win.  If Left moves on $G$ then Right can either play $G$, $i$ or $G$ and $i$.  If Right chooses to play $G$ then Left will simply respond in $G$ to force a tie, if Right plays $i$ then Left can either respond in $i$ and still tie, or play $G$ and potentially win.  If Right plays both $G$ and $i$, again Left can respond in both and tie, or play $i$ only and potentially win.  So therefore $(G\bigtriangleup i)^{SL}_F=0$.

Therefore the set of impartial games is a non-trivial monoid under the selective sum and the theorem is proven.

\end{proof}

\begin{conjecture}\label{cs1}
Not every impartial game is invertible under the selective sum.
\end{conjecture}

To prove this we have to show that for an impartial game $G$, there is no game $Y$, such that $G\triangledown Y\triangledown P\approx P$ for all impartial games $P$.  Finding such a game and proving that it has no inverse is going to be rather difficult, but nevertheless the conjecture is likely to be true.

\subsection{Sprague-Grundy Theory}

As with the other operators I will define the function in the most general possible sense.

\begin{definition}
Let $n\in O=(t_1t_2\dots t_f, p_1,\dots p_f)$ and $m\in P=(s_1s_2\dots s_e, q_1\dots q_e)$;

\begin{itemize}
\item{$\G_s(0)=0$.}
\item{$\G_s(n)=\max_{k,i}\{p_k-\G_s(n_1\plus n_2\triangledown \dots\triangledown n_{i})\}$, where $n_1+n_2+\dots +n_{i}=n-k$ and $t_k=\Sigma_{i\in S_k}2^i$.}
\item{$\G_s(n\plus m)=\max_{k,i,l,j}\{p_k - \G_s(n_1\triangledown n_2\triangledown\dots\triangledown n_i\triangledown m), q_l-\G_s(n\triangledown m_1\triangledown m_2\triangledown \dots\triangledown m_j), p_k+q_l-\G_s(n_1\triangledown n_2\triangledown \dots\triangledown n_i\triangledown m_1\triangledown m_2\triangledown\dots\triangledown m_j)\}$, where $n_1+n_2+\dots +n_i=n-k$, $t_k=\Sigma_{i\in S_k}2^i$, $m_1+m_2+\dots m_j=m-l$ and $s_l=\Sigma_{j\in R_l}2^j$.}
\end{itemize}
\end{definition}

\begin{theorem}
Suppose $O_1,\dots,O_v$ are octal games, and there are natural numbers $N_1,\dots,N_v$ such that for each$i=1,\dots,v$, $G_s(n)\geq 0$ for all $n\in O_i$ and $n\leq N_i$. Then if $n_i\in O_i$ and $n_i\leq N_i$ for each $i=1,\dots,v$, $\G_s(n_1\triangledown\dots\triangledown n_v)=\Sigma_{i=1}^v\G_s(n_i)$.
\end{theorem}

\begin{proof}  I will prove this by induction on $n_1+\dots+n_j$ for some $j$.  The base case is clearly trivial since $\G_s(0\triangledown\dots\triangledown 0)=0$ regardless of how many $0$'s there are.

So for the inductive step assume that the result holds for all $n_1+\dots+n_j\leq K$ and I will choose and $n$ and $m$ such that $n+m=K+1$, and $G_s(n)$ and $G_s(m)\geq 0$.  The reason I only choose two games $n$ and $m$ is because it makes the proof easier and it will also be clear that the same argument can be extended to any number of games.

$\G_s(n\triangledown m)=\max_{k,i,l,j}\{p_k - \G_s(n_1\triangledown\dots\triangledown n_i\triangledown m), q_l-\G_s(n\triangledown m_1\triangledown \dots\triangledown m_j), p_k+q_l-\G_s(n_1\triangledown \dots\triangledown n_i\triangledown m_1\triangledown\dots\triangledown m_j)\}$, and since $n_1+\dots n_i+m$, $m_1+\dots+m_j+n$ and $n_1\dots n_i+m_1+\dots + m_j\leq k$, then by induction, $\max_{k,i,l,j}\{p_k - \G_s(n_1\triangledown\dots\triangledown n_i\triangledown m), q_l-\G_s(n\triangledown m_1\triangledown \dots\triangledown m_j), p_k+q_l-\G_s(n_1\triangledown \dots\triangledown n_i\triangledown m_1\triangledown\dots\triangledown m_j)\}=\max\{p_k - \G_s(n_1\triangledown\dots\triangledown n_i)-\G_s(m), q_l-\G_s(n)-\G_s(m_1\triangledown \dots\triangledown m_j), p_k+q_l-\G_s(n_1\triangledown \dots\triangledown n_i) -\G_s(m_1\triangledown\dots\triangledown m_j)\}=\max\{\G_s(n)-\G_s(m), \G_s(m)-\G_s(n), \G_s(n)+\G_s(m)\}$.  

However since we know that both $\G_s(n)$ and $\G_s(m)\geq 0$, then $\max\{\G_s(n)-\G_s(m), \G_s(m)-\G_s(n), \G_s(n)+\G_s(m)\}=\G_s(n)+\G_s(m)$, as previously stated it is clear that exactly the same argument can be used for any number games and so the theorem is proven.
\end{proof}

Note that this theorem will not hold if either $\G_s(n)$ or $\G_s(m)<0$, since in that case it might be better to move on $n$ or $m$ but not both $n$ and $m$, but this is still quite a strong result and tells us quite a lot about nim variants played under the selective sum.  In the general case I make the following conjecture.

\begin{conjecture}
Let $O=(n_1n_2\dots n_t, p_1p_2\dots p_t)$ and $P=(m_1m_2\dots m_l, q_1q_2\dots q_l)$ be two finite taking-no-breaking octal games such that, there is at least one $n_s\neq 0$ or $1$, and if $n_i$ and $m_j=1,2$ or $3$ then $p_i=i$ and $q_j=j$, and $p_i=q_j=0$, otherwise, then for all $m$ there exists an $N$ such that;

$$\G_s(n+2k\triangledown m)=\G_s(n\triangledown m)$$

\noindent for all $n\geq N$ and $k$ is the largest entry in $O$ such that $n_k\neq 0,1$.
\end{conjecture}

What is interesting about this is that changing the operator does not appear to change the period, and in fact I make an even stronger conjecture;

\begin{conjecture}
Let $O=(n_1n_2\dots n_t, p_1p_2\dots p_t)$ and $P=(m_1m_2\dots m_l, q_1q_2\dots q_l)$ be two finite octal games, then if $\G_s(n\plus m)$ eventually has period $p$, $\G_s(n\triangledown m)$ also eventually has period $p$.
\end{conjecture}

So in other words what this conjecture says is that if these values are eventually periodic under the disjunctive sum, then not only are they eventually periodic under the selective sum, but they have the same period.

\section{The Sequential Join}

The sequential join was first defined by Stromquit and Ullman \cite{SU} and then studied further by Stewart \cite{FS2}.  With this operator we give all the components of a game a pre-determined order, and then play them in that order.  It is an interesting operator to look at because the structure of mis\`ere and normal play games is very similar under this operator.

\begin{definition}  The sequential join of two games $G$ and $H$ is 
defined as follows:
$$G\rhd H=\begin{cases}

&\{G^L\rhd H|G^S+H^S|G^R\rhd H\} \text{, if $G\neq \{.|G^S|.\}$}\\

&\{G^S\rhd H|G^S+H^S|G^S\rhd H\} \text{, Otherwise}

\end{cases}$$

\end{definition}

\begin{theorem}\label{seqid}
Scoring play games form a non-trivial monoid under the sequential join.
\end{theorem}

\begin{proof}  To prove this we first define a set  $I=\{i|i\rhd G\approx G\rhd i\approx G\hbox{ for all games }G\}$, and show that $I$ contains more than one element namely $\{.|0|.\}$.  So consider the game  $i=\{\{0|0|0\}|0|\{0|0|0\}\}$, as shown in the figure.
\begin{figure}[htb]
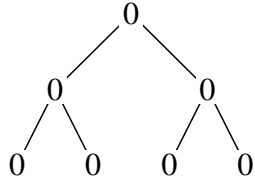

\begin{center}
\begin{graph}(3,2)

\roundnode{1}(0,0)\roundnode{2}(0.5,1)\roundnode{3}(1,0)\roundnode{4}(1.5,2)
\roundnode{5}(2,0)\roundnode{6}(2.5,1)\roundnode{7}(3,0)

\edge{1}{2}\edge{2}{3}\edge{4}{2}\edge{4}{6}\edge{5}{6}\edge{6}{7}

\nodetext{1}{0}\nodetext{2}{0}\nodetext{3}{0}\nodetext{4}{0}
\nodetext{5}{0}\nodetext{6}{0}\nodetext{7}{0}

\end{graph}
\end{center}
\caption{The game $i=\{\{0|0|0\}|0|\{0|0|0\}\}$}
\end{figure}

So first consider the game $i\rhd G$, if Left moves first on $i\rhd G$, then Right will move last on $i$, which means that Left will move first on $G$, and since the final score of $i$ is always $0$, then $(i\rhd G)^{SL}_F=G^{SL}_F$.  Similarly for the game $G\rhd i$, the players will simply play through $G$, and regardless of what happens the game $i$ cannot change the score of $G$, and therefore $(G\rhd i)^{SL}_F=G^{SL}_F$.

To show that the set is a monoid and not a group we need to demonstrate that not all games are invertible, so consider the game $Y=\{\{c|b|.\}|a|.\}\}$, and the game $G=\{e|d|f\}$.  If $Y$ is invertible this means that there exists a game $Y^{-1}$ such that $Y\rhd Y^{-1}\rhd G\approx G$ for all games $G$.  $G^{SR}_F=f$, however $(Y\rhd Y^{-1}\rhd G)^{SR}_F=a+a'+d\neq f$ and so the theorem is proven.
\end{proof}

\begin{theorem}
For any outcome classes $\mathcal{X}$, $\mathcal{Y}$ and $\mathcal{Z}$, there is a game $G\in \mathcal{X}$ and $H\in \mathcal{Y}$ such that $G\rhd H\in \mathcal{Z}$.
\end{theorem}

\begin{proof}  To prove this let $G=\{\{c|b|.\}|a|.\}$ and $H=\{H^L|d|H^R\}$, where $H^L$ and $H^R\neq \emptyset$, then $G^{SL}_F=a$, $G^{SR}_F=b$, $(G\rhd H)^{SL}_F=a+d$ and $(G\rhd H)^{SR}_F=b+d$.  Since $d$ is not dependent on $H^{SL}_F$ and $H^{SR}_F$, and can be any real number, then we can pick $a$, $b$ and $d$, so that $G$ and $H$ are in any outcome class and $G\rhd H$ is any outcome class.  Therefore the theorem is proven.  
\end{proof}

\subsection{Impartial Games}

\begin{theorem}
Impartial games for a non-trivial monoid under the sequential join.
\end{theorem}

\begin{proof}  From the proof of theorem \ref{seqid} we know that there is a non-trivial identity set, so to prove this we simply need to show that there is a game $G$ that is not invertible.  So consider the game $G=\{1,\{0|0|0\}|0|\{0|0|0\},-1\}$.  Let $Y$ be the inverse of $G$, then this implies that $G\rhd Y\rhd P\approx P$ for all impartial games $P$.

So let $P=\{.|0|.\}$, and consider the game $G\rhd Y\rhd P$.  If Left moves first and moves to the game $1\rhd Y\rhd P$, then his implies that $-1$ is one of the Right options of $G$, since if Right moves to $-1$ on $Y$ then Left will move first on $P$ and $G\rhd Y$ will not change the final score of $P$.  But $Y$ is impartial, so this implies that $1$ is a Left option of $Y$.  So therefore if Left moves to the game $\{0|0|0\}\rhd Y\rhd P$, then this means that Right must move to the game $0\rhd Y\rhd P$, and Left will move first on $Y$, and Left can choose the option $1$ and hence win $G\rhd Y\rhd P$, i.e. $G\rhd Y\rhd P\not \approx P$ which is a contradiction.

So this means that $G$ is not invertible, and therefore the set of impartial games form a non-trivial monoid under the sequential join and the theorem is proven.
\end{proof}

\subsection{Sprague-Grundy Theory}

When consider the sequential join it doesn't really make sense to look at taking and breaking games, because once you break the heap into two or more smaller heaps we have to define the order that we play the two new heaps in.  Since this is a rather difficult issue to resolve I will not be considering it in this paper.

\begin{definition}
Let $n\in O=(t_1t_2\dots t_f, p_1,\dots p_f)$ and $m\in P=(s_1s_2\dots s_e, q_1\dots q_e)$, be two taking no breaking games;

\begin{itemize}
\item{$\G_s(0)=0$.}
\item{$\G_s(n\rhd m)=\begin{cases}
&\max\{p_k-\G_s(n-k\rhd m)\} \text{, if $n\neq 0$}\\
&\max\{q_l-\G_s(n\rhd m-l)\} \text{, Otherwise}
\end{cases}$}.
\end{itemize}
\end{definition}

There is not really a lot to say about this operator, other than to make the following conjecture;

\begin{conjecture}
Let $O=(n_1n_2\dots n_t, p_1p_2\dots p_t)$ and $P=(m_1m_2\dots m_l, q_1q_2\dots q_l)$ be two finite octal games, then if $\G_s(n\plus m)$ eventually has period $p$, $\G_s(n\rhd m)$ also eventually has period $p$.
\end{conjecture}

This conjecture seems quite a reasonable one due to the nature of the operator.  By playing the heaps in order, it means that $m$ cannot change the period of $n$.  However since it is very hard to even prove that $\G_s(n+p)=\G_s(n)$ for all $n$ large enough, a proof of this conjecture will also be very difficult.

\section{Conclusion}

In this paper I really only examined the basic structure of the three operators that I looked at.  There are of course plenty of other questions that I could have looked at.  I feel that the most interesting thing was looking at the function $\G_s(n)$ under each of the different operators as it appears that the function settles down into the same period regardless of the operator being used.

I hope to be able to prove all of the conjectures made in this paper, and I feel that a proof of them would tell us a lot about the function $\G_s(n)$ and the nature of octal games under scoring play rules.

\end{document}